\documentclass[a4paper,12pt]{article}
\usepackage{amsmath,amssymb,amsthm}
\usepackage{graphicx,psfrag}

\numberwithin{equation}{section}

 \DeclareMathOperator{\id}{id}
 \DeclareMathOperator{\im}{Im}
 \DeclareMathOperator{\tr}{tr}
 \DeclareMathOperator{\vol}{vol}
 \DeclareMathOperator{\Gal}{Gal}

 \newcommand{\C}{\mathbb{C}}
 \newcommand{\N}{\mathbb{N}}
 \newcommand{\Q}{\mathbb{Q}}
 \newcommand{\R}{\mathbb{R}}
 \newcommand{\Z}{\mathbb{Z}}
 \newcommand{\pii}{\pi\sqrt{-1}}
 \newcommand{\piN}[1]{\frac{#1\pi}{N}}
 \newcommand{\Npi}{\frac{N}{\pi}}
 \newcommand{\hN}{\frac{N}2}
 \newcommand{\qN}{\frac{N}4}
 \newcommand{\lognorm}[1]{\log \Big| #1 \Big|}
 \newcommand{\rt}{e^{\frac{2\pii}N}}
 \newcommand{\J}[1]{J_{#1,N}(\rt)}
 \newcommand{\hJ}[2]{\hat{J}_{#1,#2}(\rt)}
 \newcommand{\logJ}[1]{\log \big| \J{#1} \big|}
 \newcommand{\cv}[1]{\vol(S^3 \setminus #1)}
 \newcommand{\sv}[1]{\| S^3 \setminus #1 \|}

\newtheorem{thm}{Theorem}[section]
\newtheorem{lem}[thm]{Lemma}
\newtheorem{prop}[thm]{Proposition}
\newtheorem{conj}[thm]{Conjecture}
\theoremstyle{remark}
\newtheorem{rem}[thm]{Remark}

\begin{document}

\title{Proof of the volume conjecture for Whitehead doubles of a family of torus knots}
\author{Hao Zheng \\
    {\small Department of Mathematics, Zhongshan University} \\
    {\small Guangzhou 510275, China} \\
    {\small E-mail: zhenghao@mail.sysu.edu.cn}
}
\date{}
\maketitle

\begin{abstract}
A technique to calculate the colored Jones polynomials of
satellite knots, illustrated by the Whitehead doubles of knots, is
presented. Then we prove the volume conjecture for Whitehead
doubles of a family of torus knots and show some interesting
observations.
\end{abstract}

\begin{small}

\smallskip

\end{small}

\section{Introduction}

The volume conjecture was proposed by Kashaev and reformulated and
refined by Murakami and Murakami as follows.

\begin{conj}[Kashaev \cite{Kashaev}, Murakami-Murakami \cite{MM}]
For any knot $K$,
\begin{equation}\label{eqn:conj}
  2\pi \lim_{N\to\infty} \frac {\logJ{K}} {N} = v_3 \sv{K}
\end{equation}
where $J_{K,N}$ is the (normalized) colored Jones polynomial of
$K$, $\sv{K}$ is the simplicial volume of the complement of $K$
and $v_3$ is the volume of the ideal regular tetrahedron.
\end{conj}

Recall that $v_3\sv{K}$ is nothing but the sum of the hyperbolic
volumes of hyperbolic pieces in the JSJ-decomposition of the
complement of $K$. In Kashaev's original form, the knot $K$ is
hyperbolic and the equation is in terms of the quantum dilogarithm
invariant and the hyperbolic volume of the complement of $K$.

The conjecture is marvellous in the sense that it reveals the
topological meaning of the quantum invariants of knots which is
quite unobvious from definition. However, it also turns out to be
rather hard to be proved. Till now, besides positive numerical
evidences (ref. \cite{Hikami,MMOTY}) for some hyperbolic knots,
only the cases of torus knots (Kashaev-Tirkkonen \cite{KT}) and
the simplest hyperbolic knot, the figure 8 knot (ref.
\cite{Kashaev}) have been verified.

In view of the compatible behavior of both sides of the
conjectured equation (\ref{eqn:conj}) under connect sum
\begin{eqnarray}
 && J_{K_1 \sharp K_2,N} = J_{K_1,N} \cdot J_{K_2,N}, \\
 && \sv{K_1 \sharp K_2} = \sv{K_1} + \sv{K_2},
\end{eqnarray}
the volume conjecture, in fact, may be reduced to the
consideration of prime knots. By Thurston's Hyperbolization
Theorem (ref. \cite{Thurston}), the prime knots further fall into
three families: torus knots, hyperbolic knots and satellite knots.

In this article, we deal with the conjecture by examining a
special case of the third family, the Whitehead doubles of torus
knots. The approach is emphasized on the relation between the
colored Jones polynomial of a satellite knot and those of the
associated companion knot and pattern link. In particular, we show
a technique to calculate the colored Jones polynomial of satellite
knots by cutting and gluing method.

\begin{center}
    \includegraphics[scale=.7]{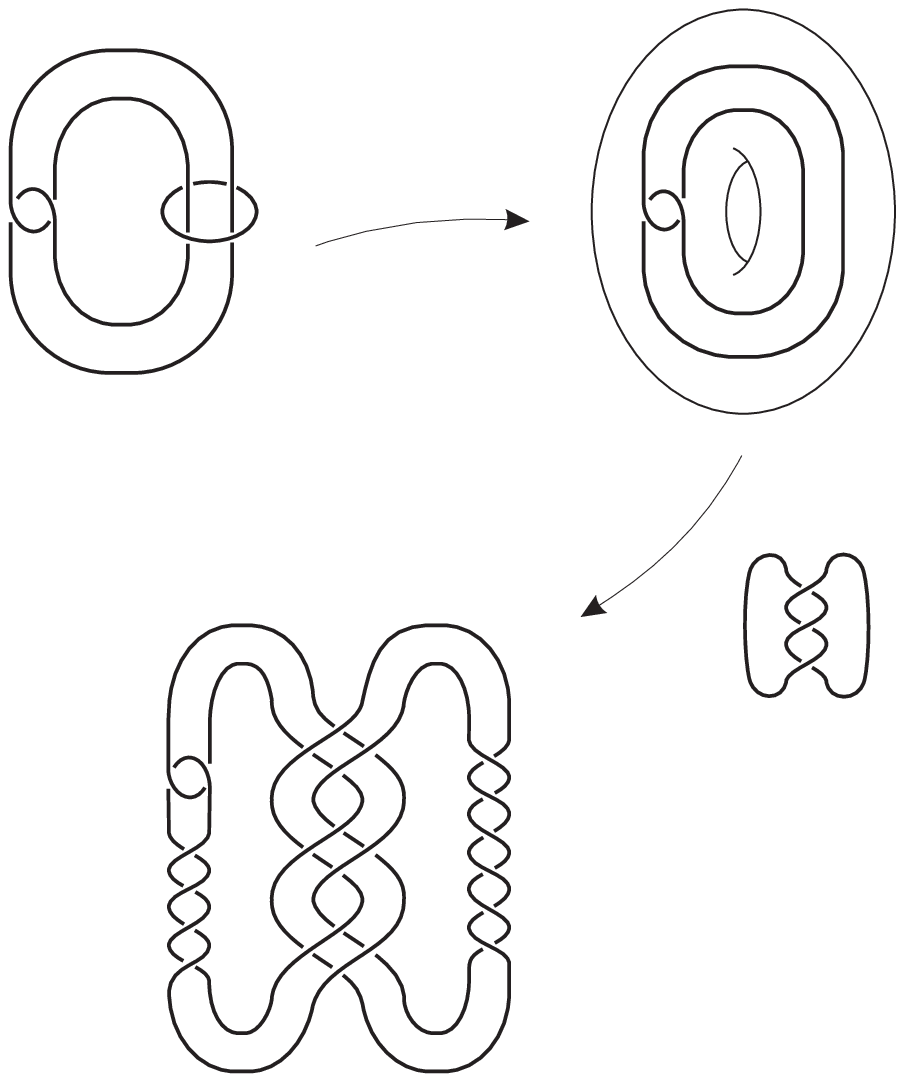}
\end{center}

A Whitehead double of a knot $K$ is a knot obtained as follows.
Remove the regular neighborhood of one component of the Whitehead
link from $S^3$ thus get a knot inside a torus, then knot the
torus in the shape of a knot $K$.

Note that, when $K$ is nontrivial, a Whitehead double $K'$ of $K$
is a satellite knot whose complement contains an obvious essential
torus $T^2$. Cutting along the torus, we get
\begin{equation}
  (S^3 \setminus K') \setminus T^2 \cong
  (S^3 \setminus \text{Whitehead link}) \cup (S^3 \setminus K)
\end{equation}
thus
\begin{equation}
  \sv{K'} = \sv{\text{Whitehead link}} + \sv{K}.
\end{equation}
In particular, if $K$ is a nontrivial torus knot, the complement
of $K$ is Seifert fibred and the complement of the Whitehead link
is hyperbolic, hence
\begin{equation}
  v_3\sv{K'} = \cv{\text{Whitehead link}}.
\end{equation}

The article proceeds as follows. First, we calculate the colored
Jones polynomials of the twisted Whitehead links and the Whitehead
doubles of knots in section \ref{sec:tangle}. Next, as a
warming-up we prove in the next two consecutive sections the
following two theorems, of which the former one is, in fact, the
volume conjecture for twisted Whitehead links and both extends the
estimation (\ref{eqn:conj}) to the second order.

\begin{thm}\label{thm:wl}
For every twisted Whitehead link $L$, we have
\begin{equation}
  2\pi \logJ{L} = \cv{L} \cdot N  + 3\pi \log N + O(1)
\end{equation}
as $N\to\infty$.
\end{thm}

\begin{thm}\label{thm:torus}
For every nontrivial torus knot $T(p,q)$ with $q=2$, we have
\begin{equation}\label{eqn:torus}
  2\pi \logJ{T(p,q)} = 3\pi \log N + O(1)
\end{equation}
as $N\to\infty$.
\end{thm}

Then we prove the main theorem in Section \ref{sec:wd} and show
some observations in the final section.
\begin{thm}\label{thm:wd}
If $K$ is a Whitehead double of a nontrivial torus knot $T(p,q)$
with $q=2$, then
\begin{equation}\label{eqn:main}
  2\pi \logJ{K} = v_3\sv{K} \cdot N  + 4\pi \log N + O(1)
\end{equation}
as $N\to\infty$. In particular, the volume conjecture is true for
$K$.
\end{thm}

\medskip

\begin{rem}
In their proof of the volume conjecture for torus knots, Kashaev
and Tirkkonen \cite{KT} derived the following estimation
\begin{equation}
  2\pi \logJ{T(p,q)} = O(\log N).
\end{equation}
But improving the estimation to (\ref{eqn:torus}) requires a
nonvanishing proposition on number theory (see Proposition
\ref{prop:novanish}) to which both Theorem \ref{thm:torus} and
Theorem \ref{thm:wd} are reduced in this article. With a technical
condition $q=2$ we proved the nonvanishing proposition in section
\ref{sec:torus}. A complete proof has been beyond the scope of the
article. We only mention here that our technique can be sharpened
to prove the nonvanishing proposition, hence both theorems, at
least for the cases that both $p,q$ are odd or one of them is a
power of $2$.
\end{rem}

\begin{rem}
It is noteworthy that the coefficient ``$4\pi$" of the second term
in the asymptotic expansion (\ref{eqn:main}) disagrees with the
observation due to Hikami \cite{Hikami}
\begin{equation}
  2\pi \logJ{K} = v_3\sv{K} \cdot N  + 3\pi \log N + O(1)
\end{equation}
for many prime knots $K$.
\end{rem}

\section{Calculation of colored Jones polynomial}
\label{sec:tangle}

In this section, we calculate the colored Jones polynomials of the
twisted Whitehead link $WL(r)$ and the Whitehead double $WD(K,r)$
of a knot $K$.
\begin{center}
    \psfrag{a}[c]{WL($r$)}
    \psfrag{b}[c]{WD($K$,$r$)}
    \psfrag{c}[c]{double($K$)}
    \psfrag{d}[c]{$r$ twists}
    \includegraphics[scale=.8]{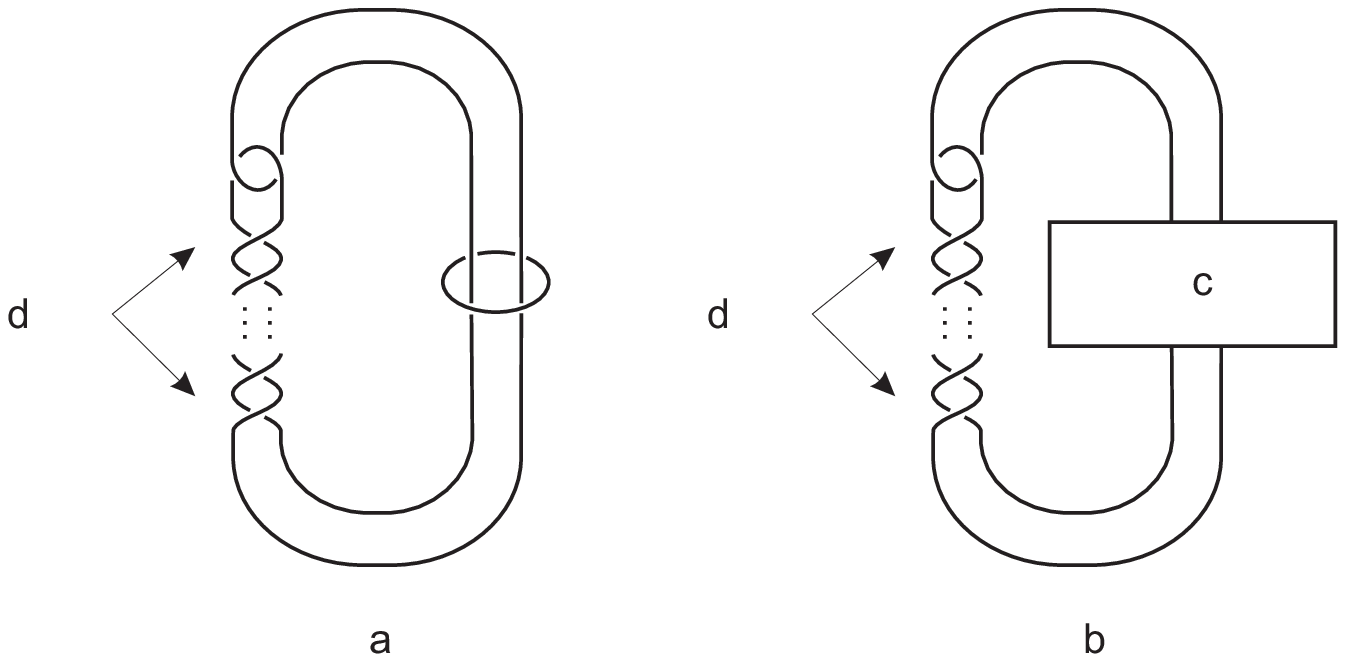}
\end{center}
In the figure, $double(K)$ denotes the (2,2)-tangle obtained by
doubling the knot $K$ to a link with zero linking number and then
removing a pair of parallel segments.

Our trick is cutting the link diagrams into (2,2)-tangles and
gluing the tangle invariants together.

\begin{center}
    \psfrag{a}[c]{twist}
    \psfrag{b}[c]{belt}
    \psfrag{c}[c]{double($K$)}
    \psfrag{d}[c]{clasp}
    \includegraphics[scale=1]{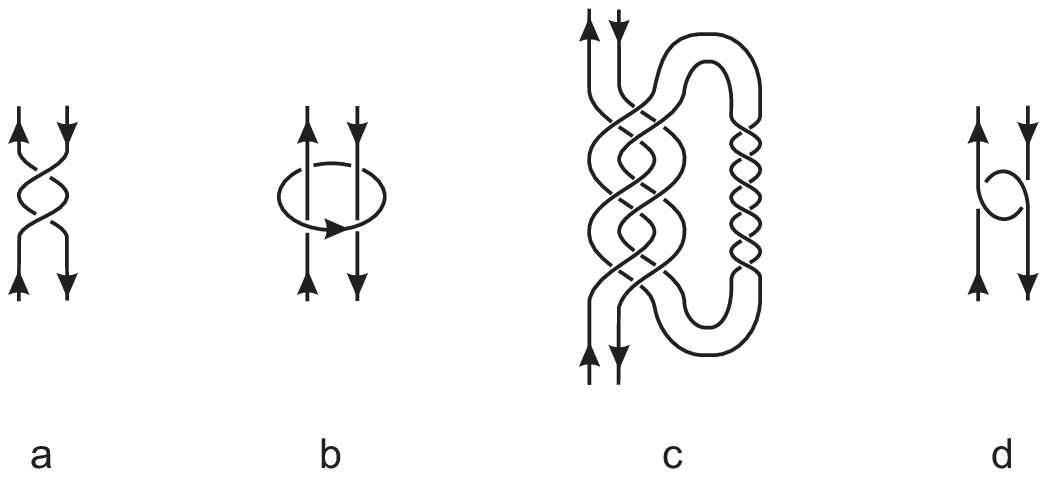}
\end{center}

Colored Jones polynomial is also defined for tangles, but, instead
of a Laurent polynomial of $t$, it is in general a module
homomorphism of $U_q(sl_2)$ (choose $t = q^2$). Especially, the
colored Jones polynomial of a (2,2)-tangle is a module
homomorphism
\begin{equation}
  V_N \otimes V_N \to V_N \otimes V_N
\end{equation}
where $V_N$ is the $N$ dimensional irreducible representation of
$U_q(sl_2)$.

Note that the tensor product admits the decomposition
\begin{equation}
  V_N \otimes V_N = \bigoplus_{n=0}^{N-1} V_{2n+1}.
\end{equation}
A straightforward calculation shows that the (framing independent,
unnormalized) colored Jones polynomials of the tangles are

\begin{eqnarray}
  && \tilde{J}_{twist,N} = \bigoplus_{n=0}^{N-1} t^{n(n+1)}
  \cdot \id_{V_{2n+1}}, \\
  && \tilde{J}_{belt,N} = \bigoplus_{n=0}^{N-1}
  \frac {t^{N(2n+1)/2}-t^{-N(2n+1)/2}} {t^{(2n+1)/2}-t^{-(2n+1)/2}}
  \cdot \id_{V_{2n+1}}, \\
  && \tilde{J}_{double(K),N} = \bigoplus_{n=0}^{N-1} J_{K,2n+1}
  \cdot \id_{V_{2n+1}}, \\
  && \tilde{J}_{clasp,N} = \bigoplus_{n=0}^{N-1} \xi_{N,n}
  \cdot \id_{V_{2n+1}},
\end{eqnarray}
where
\begin{equation}
  \xi_{N,n} = t^{(N^2-1)/2 + N(N-1)/2} \sum_{i=0}^{N-1-n} t^{-N(i+n)}
  \prod_{j=1}^{n} \frac {(1-t^{N-i-j})(1-t^{i+j})} {1-t^j}.
\end{equation}

Combining the tangle invariants together, one has
\begin{equation}
  J_{WL(r),N}
  = \sum_{n=0}^{N-1}
  \frac {t^{(2n+1)/2}-t^{-(2n+1)/2}} {t^{N/2}-t^{-N/2}} \cdot
  t^{rn(n+1)} \cdot \xi_{N,n} \cdot
  \frac {t^{N(2n+1)/2}-t^{-N(2n+1)/2}} {t^{(2n+1)/2}-t^{-(2n+1)/2}},
\end{equation}
and
\begin{equation}\label{eqn:Jwd}
  J_{WD(K,r),N}
  = \sum_{n=0}^{N-1}
  \frac {t^{(2n+1)/2}-t^{-(2n+1)/2}} {t^{N/2}-t^{-N/2}} \cdot
  t^{rn(n+1)} \cdot \xi_{N,n} \cdot J_{K,2n+1}.
\end{equation}
Note that, in the expression of $J_{WD(K,r),N}$, the factor
$J_{K,2n+1}$ is contributed by the companion knot $K$ and the
other part is precisely obtained from the expression of
$J_{WL(r),N}$ by removing the factor contributed by the belt
tangle.

\section{Proof of Theorem \ref{thm:wl}}
\label{sec:wl}

Let $L$ denote the twisted Whitehead link $WL(r)$. Setting
$t=\rt$, we have
\begin{equation}\label{eqn:wl}
\begin{split}
  \J{L}
  & = -t^{-1/2}
  \sum_{n=0}^{N-1} (2n+1) t^{rn(n+1)} \sum_{i=0}^{N-1-n}
  \prod_{j=1}^{n} \frac {(1-t^{-i-j})(1-t^{i+j})} {1-t^j} \\
  & = -e^{-\frac{\pii}{N}}
  \sum_{n=0}^{N-1} (2n+1) a_n^{4r-1} \sum_{i=0}^{N-1-n} S_{n,i}
\end{split}
\end{equation}
where
\begin{equation}
  a_n = e^{\frac{n(n+1)}{2N}\pii - \frac{n}{2}\pii}
  = e^{\frac{n(n+1-N)}{2N}\pii}
\end{equation}
and
\begin{equation}
   S_{n,i} = \prod_{j=1}^{n}
  \frac {4\sin^2\piN{(i+j)}} {2\sin\piN{j}}.
\end{equation}

First, we prepare a lemma to estimate the norm factor $S_{n,i}$.
Put
\begin{equation}
  s_n = - \sum_{j=1}^n \lognorm{2\sin\piN{j}}
\end{equation}
and let
\begin{equation}
  L(x) = - \int_0^x \log|2\sin u| du
\end{equation}
be the Lobachevsky function.

\begin{lem}\label{lem:sest}
For $0 < \alpha < 1$ we have uniform estimations
\begin{equation}\label{eqn:sest1}
  s_m - s_n = \Npi L(\piN{m}) - \Npi L(\piN{n}) + O(N^{-1})(m-n)
\end{equation}
on $\frac\alpha2 N < n < m < (1-\frac\alpha2)N$,
\begin{equation}\label{eqn:sest2}
  s_n = \Npi L(\piN{n}) - \frac12 \log n + O(1)
\end{equation}
on $0 < n < \alpha N$ and
\begin{equation}\label{eqn:sest3}
  s_n = \Npi L(\piN{n}) - \frac12 \log (N - n) + O(1)
\end{equation}
on $(1-\alpha)N < n < N$.
\end{lem}

\begin{proof}
We have
\begin{equation}
\begin{split}
  & - \lognorm{2\sin\piN{j}} +
  \Npi \int_\piN{(j-1)}^\piN{j} \log|2\sin u| du \\
  = & - \lognorm{2\sin\piN{j}} +
  \Npi \int_0^\piN{} \lognorm{2\sin(\piN{j}-u)} du \\
  = & \Npi \int_0^\piN{}
  \lognorm{ \frac {\sin(\piN{j}-u)} {\sin\piN{j}} } du.
\end{split}
\end{equation}

Since
\begin{equation}
  \lognorm{\frac {\sin(x-u)} {\sin x}} = O(u)
\end{equation}
uniformly on $x \in [\frac\alpha2\pi,(1-\frac\alpha2)\pi]$ as $u
\to 0$, the first estimation follows as
\begin{equation}
\begin{split}
  & \quad s_m -s_n - \Npi L(\piN{m}) + \Npi L(\piN{n}) \\
  & = \sum_{j=n+1}^m \Npi \int_0^\piN{}
  \lognorm{ \frac {\sin(\piN{j}-u)} {\sin\piN{j}} } du \\
  & = O(N^{-1}) (m-n).
\end{split}
\end{equation}

Note that
\begin{equation}
  \frac {\sin(x-u)} {x-u} \cdot \frac {x} {\sin x}
  = 1 + O(u)
\end{equation}
thus
\begin{equation}
  \lognorm{\frac {\sin(x-u)} {\sin x}}
  = \lognorm{\frac {x-u} {x}} + O(u)
\end{equation}
uniformly on $x \in [-\alpha\pi,\alpha\pi] \setminus \{0,u\}$ as
$u \to 0$. It follows that
\begin{equation}
\begin{split}
  s_n - \Npi L(\piN{n})
  & = \sum_{j=1}^n \Npi \int_0^\piN{}
  \lognorm{ \frac {\sin(\piN{j}-u)} {\sin\piN{j}} } du \\
  & = \sum_{j=1}^n \Npi \int_0^\piN{}
  \lognorm{ \frac {\piN{j}-u} {\piN{j}} } du
  + n O(N^{-1}) \\
  & = \sum_{j=1}^n \Big( - (j-1) \log \frac{j-1}{j} - 1 \Big)
  + O(1) \\
  & = - \log \frac{n!}{n^n} - n + O(1)
\end{split}
\end{equation}
uniformly on $0 < n < \alpha N$. Thanks to Sterling series
\begin{equation}
  \log n! = n \log n - n + \frac12 \log n + \frac12 \log2\pi + \cdots,
\end{equation}
the second estimation holds.

To see the third estimation, one notices that
\begin{equation}
  L(x) + L(\pi-x) = 0
\end{equation}
and
\begin{equation}
  s_{n-1} + s_{N-n} = s_{N-1}.
\end{equation}
In particular, we have
\begin{equation}
  L(\frac\pi2) = 0
\end{equation}
and, by the second estimation,
\begin{equation}
  s_{N-1} = s_{[\frac{N-1}2]} + s_{[\hN]}
  = - \log N + O(1).
\end{equation}
Therefore,
\begin{equation}
\begin{split}
  s_n & = s_{N-1} - s_{N-n} - \lognorm{2\sin\piN{n}} \\
  & = - \log N + \Npi L(\piN{n}) + \frac12 \log (N - n)
  - \log \piN{2(N-n)} + O(1) \\
  & = \Npi L(\piN{n}) - \frac12 \log (N - n) + O(1)
\end{split}
\end{equation}
uniformly on $(1-\alpha)N < n < N$.
\end{proof}

From the second and the third estimations of above lemma, we have
\begin{equation}\label{eqn:Sest}
  \log S_{n,i}
  = -2s_{n+i} + 2s_i + s_n
  = \Npi f \big( \piN{n},\piN{i} \big) + O(\log N)
\end{equation}
uniformly on $0 \leq n,i,n+i < N$, where
\begin{equation}
  f(x,y) = -2L(x+y)+2L(y)+L(x).
\end{equation}
The function $f(x,y)$ has a unique critical point
$(\frac\pi2,\frac\pi4)$ in the region $0 \leq x,y,x+y \leq \pi$,
at which $f$ reaches maximum
\begin{equation}\label{eqn:fexp}
  f \big( \frac\pi2,\frac\pi4 \big)
  = 4L \big( \frac\pi4 \big)
\end{equation}
and expands as
\begin{equation}\label{eqn:fexp}
  f \big( x+\frac\pi2,y+\frac\pi4 \big)
  = f \big( \frac\pi2,\frac\pi4 \big) - (x^2+2xy+2y^2) + \cdots.
\end{equation}

Notice that the phase factor $a_n$ is also steady near $\hN$. In
what follows, the summation (\ref{eqn:wl}) is expected to be
dominated by the summands whose index $(n,i)$ is near $(\hN,\qN)$.
Indeed, this is the case as demonstrated by the next pair of
lemmas.

\begin{lem}\label{lem:main1}
For any $\frac12 < \delta < 1$ there exist $\epsilon>0$ and $C>0$
such that
\begin{equation}
  S_{n,i} < C e^{-\epsilon N^{2\delta-1}} S_{[\hN],[\qN]}
\end{equation}
for $|n-\hN|+|i-\qN| \geq N^\delta$
\end{lem}

\begin{proof}
Since $f$ has a unique critical point $(\frac\pi2,\frac\pi4)$ in
the region $0 \leq x,y,x+y \leq \pi$, we have
\begin{equation}
  f(x,y)
  \leq \max_{|x'-\frac\pi2|+|y'-\frac\pi4| = \pi N^{\delta-1}}
  f(x',y')
\end{equation}
for $|x-\frac\pi2|+|y-\frac\pi4| \geq \pi N^{\delta-1}$. By
(\ref{eqn:fexp}), there exist $\epsilon>0$ and $C'>0$ such that
\begin{equation}
  \max_{|x'-\frac\pi2|+|y'-\frac\pi4| = \pi N^{\delta-1}}
  f(x',y')
  < f(\frac\pi2,\frac\pi4) - 2\pi\epsilon (N^{\delta-1})^2 + C'
\end{equation}
Therefore, by (\ref{eqn:Sest}) there exists $C''>0$ such that
\begin{equation}
  \log S_{n,i}
  < \log S_{[\hN],[\qN]} - \epsilon N^{2\delta-1} + C''
\end{equation}
for $|n-\hN|+|i-\qN| \geq N^\delta$.
\end{proof}

\begin{lem}\label{lem:main2}
For any $\alpha \geq 0$, $\beta\in\R$ and $\frac12 < \delta <
\frac23$ there exists a nonzero constant $C \in \C$ such that
\begin{equation}
  \sum_{|n-\hN|+|i-\qN| < N^\delta}
  (2n+1)^\alpha a_n^\beta S_{n,i}
  = C N^{\alpha+1} e^{-\frac{\beta N}8\pii}
  S_{[\hN],[\qN]} (1+O(N^{3\delta-2})).
\end{equation}
\end{lem}
\begin{proof}
For simplicity, we use the notation $n'=n-\hN$, $i'=i-\qN$ in the
proof. Note that
\begin{equation}\label{eqn:gauss}
\begin{split}
  & \sum_{|n-\hN|+|i-\qN| < N^\delta}
  (2n+1)^\alpha
  e^{-\piN{}(n'^2+2n'i'+2i'^2)+\frac{\beta n'^2}{2N}\pii} \\
  & = \int_{|x|+|y|<N^{\delta-\frac12}}
  N^{\alpha+1}
  e^{-\pi(x^2+2xy+2y^2)+\frac{\beta x^2}2\pii} dxdy
  (1+O(N^{\delta-1})) \\
  & = \int_{\R^2} N^{\alpha+1}
  e^{-\pi(x^2+2xy+2y^2)+\frac{\beta x^2}2\pii} dxdy
  (1+O(N^{\delta-1})).
\end{split}
\end{equation}
By (\ref{eqn:sest1}) and (\ref{eqn:fexp}) we have
\begin{equation}
\begin{split}
  \log S_{n,i} - \log S_{[\hN],[\qN]}
  & = \Npi f(\piN{n},\piN{i}) - \Npi f(\frac\pi2,\frac\pi4)
  + O(N^{\delta-1}) \\
  & = - \piN{} (n'^2+2n'i'+2i'^2)
  +  O(N^{3\delta-2})
\end{split}
\end{equation}
uniformly on $|n-\hN|+|i-\qN| < N^\delta$. Moreover, on the same
region we have the uniform estimation
\begin{equation}
\begin{split}
  & a_n = e^{ \frac{n(n+1-N)}{2N}\pii }
  = e^{ (\frac{n'^2}{2N} - \frac{N}8 + \frac{n}{2N}) \pii } \\
  & \quad\quad\quad\quad\quad
  = e^{ (\frac{n'^2}{2N} - \frac{N}8 + \frac14) \pii }
  (1+O(N^{\delta-1})).
\end{split}
\end{equation}
Therefore, by (\ref{eqn:gauss}),
\begin{equation}
\begin{split}
  & \sum_{|n-\hN|+|i-\qN| < N^\delta}
  (2n+1)^\alpha a_n^\beta S_{n,i}
  = \int_{\R^2}
  e^{-\pi(x^2+2xy+2y^2)+\frac{\beta x^2}2\pii} dxdy \\
  & \quad\quad\quad\quad\quad\quad\quad\quad\quad\quad
  \cdot N^{\alpha+1} e^{\beta(-\frac{N}8+\frac14)\pii}
  S_{[\hN],[\qN]} (1+O(N^{3\delta-2})).
\end{split}
\end{equation}
To conclude the lemma it suffices to choose
\begin{equation}
  C = e^{\frac\beta4\pii}
  \int_{\R^2} e^{-\pi(x^2+2xy+2y^2)+\frac{\beta x^2}2\pii} dxdy.
\end{equation}
\end{proof}

It follows from Lemma \ref{lem:main1} and Lemma \ref{lem:main2}
that, in the same notations as the lemmas,
\begin{equation}
\begin{split}
  & \sum_{|n-\hN|+|i-\qN| \geq N^\delta} (2n+1) a_n^{4r-1} S_{n,i}
  = N^3 e^{-\epsilon N^{2\delta-1}} S_{[\hN],[\qN]} O(1), \\
  & \sum_{|n-\hN|+|i-\qN| < N^\delta} (2n+1) a_n^{4r-1} S_{n,i}
  = N^2 S_{[\hN],[\qN]} e^{O(1)},
\end{split}
\end{equation}
so
\begin{equation}
  \logJ{L} = \log (N^2 S_{[\hN],[\qN]}) + O(1).
\end{equation}
From (\ref{eqn:sest2}) we also have
\begin{equation}
  \log S_{[\hN],[\qN]}
  = \frac{4N}{\pi} L(\frac{\pi}{4}) - \frac12 \log N + O(1).
\end{equation}
Therefore,
\begin{equation}
\begin{split}
  2\pi \logJ{L}
  & = 8L(\frac{\pi}{4}) \cdot N + 3\pi \log N + O(1) \\
  & = \cv{L} \cdot N + 3\pi \log N + O(1)
\end{split}
\end{equation}
as $N\to\infty$. In the last row, we used the fact
\begin{equation}
  \cv{L} = \cv{\text{Whitehead link}}
  = 8L(\frac{\pi}{4}).
\end{equation}

\section{Proof of Theorem \ref{thm:torus}} \label{sec:torus}

The colored Jones polynomial of the torus knot $T(p,q)$ was
calculated in \cite{Morton} as
\begin{equation}
  J_{T(p,q),n}
  = \frac {t^{-pq(n^2-1)/4}} {t^{n/2}-t^{-n/2}}
  \sum_{k=-(n-1)/2}^{(n-1)/2} t^{pk(qk+1)} (t^{qk+1/2}-t^{-qk-1/2}).
\end{equation}
We put
\begin{equation}
  A_{p,q}^\pm(N,k)
  = \sum_{j=1}^{pq-1} (\pm1)^j e^{-\frac{Nj^2}{2pq}\pi\sqrt{-1}}
  j^{2k} \sin\frac{j\pi}{p} \sin\frac{j\pi}{q}.
\end{equation}
Note that
\begin{equation}
  A_{p,q}^\pm(N,k) = A_{p,q}^\mp(N+2pq,k) = A_{p,q}^\pm(N+4pq,k),
\end{equation}
so
\begin{equation}\label{eqn:A}
\begin{split}
  A_{p,q}^{(-1)^n}(N,k) = A_{p,q}^+(N+2npq,k), \\
  A_{p,q}^+(N,k) = A_{p,q}^{(-1)^n}(N+2npq,k).
\end{split}
\end{equation}

In \cite{KT}, an estimation of $\J{T(p,q)}$ was derived as
\begin{equation}\label{eqn:torus_estN}
  \J{T(p,q)}
  = 2 e^{ -pq\frac{N^2-1}{2N}\pii }  \frac{N^{3/2}}{(2pq)^{3/2}}
  e^{ -(\frac{p}{q}+\frac{q}{p}) \frac{\pii}{2N} + \frac\pii4 }
  A_{p,q}^{(-1)^{N-1}}(N,1) + O(1).
\end{equation}
In view of the periodicity of $A_{p,q}^\pm$, to establish the
theorem it suffices to show that $A_{p,q}^{(-1)^{N-1}}(N,1)$ never
vanishes if $q=2$, or equivalently by (\ref{eqn:A}),

\begin{prop}\label{prop:novanish}
Let $p,q \geq 2$ be coprime integers with $q=2$. Then for every
integer $N$,
\begin{equation}
  A_{p,q}^+(N,1)
  = \sum_{j=1}^{pq-1} e^{-\frac{Nj^2}{2pq}\pi\sqrt{-1}}
  j^2 \sin\frac{j\pi}{p} \sin\frac{j\pi}{q} \neq 0.
\end{equation}
\end{prop}

The proof of the proposition is purely arguments on elementary
algebraic number theory. In the following, we write $\zeta_n =
e^{\frac{2\pii}{n}}$ for each $n\in\N$. An algebraic number field
means a finite extension of $\Q$ contained in $\C$.

For any finite extension $E/K$ of field, one has a $K$-linear map
$\tr_{E/K} : E \to K$, called the {\em trace function}, which
values on $x \in E$ the trace of the $K$-linear transformation
$\rho_x : E \to E$ given by $\rho_x(z) = xz$.

\begin{lem}\label{lem:A1}
Let $\alpha$ be a prime, $k,l \in \N$ and $K$ be an algebraic
number field such that $K \cap \Q(\zeta_{\alpha^{k+l}}) = \Q$.
Then we have
\begin{equation}
  \tr_{K(\zeta_{\alpha^{k+l}})/K(\zeta_{\alpha^l})}
  (\zeta_{\alpha^{k+l}}^n) =
  \left\{ \begin{array}{ll}
    0, & \alpha^k \nmid n, \\
    \alpha^k \cdot \zeta_{\alpha^{k+l}}^n, & \alpha^k \mid n.
  \end{array} \right.
\end{equation}
\end{lem}

\begin{proof}
The field extension $K(\zeta_{\alpha^{k+l}})/K(\zeta_{\alpha^l})$
has a basis $\{ \zeta_{\alpha^{k+l}}^i \mid 0 \leq i < \alpha^k
\}$ on which the diagonal of the matrix of
$\zeta_{\alpha^{k+l}}^n$ consists of only $0$ if $\alpha^k \nmid
n$, or $\zeta_{\alpha^{k+l}}^n$ otherwise.
\end{proof}

\begin{lem}\label{lem:A2}
Let $\alpha$ be a prime and $K$ be an algebraic number field such
that $K \cap \Q(\zeta_{\alpha}) = \Q$. Then
\begin{equation}
  \sum_{j=0}^{\alpha-1} c_j \cdot \zeta_{\alpha}^j = 0,
\end{equation}
for $c_j \in K$ if and only if the $c_j$'s are identical.
\end{lem}

\begin{proof}
On one hand, the field extension $K(\zeta_{\alpha})/K$ has a basis
$\{ 1, \zeta_{\alpha}, \zeta_{\alpha}^2, \dots,
\zeta_{\alpha}^{\alpha-2} \}$. On the other hand, we have
$\sum_{j=0}^{\alpha-1} \zeta_\alpha^j = 0$. Therefore, the
summation vanishes if and only if the $c_j$'s are identical to
$c_{\alpha-1}$.
\end{proof}

Thanks to the next lemma, we are able to eliminate the Gaussian
exponential appearing in the expression of $A_{p,q}^\pm$.

\begin{lem}\label{lem:A4}
Let $\alpha$ be an odd prime, $l \in \N$ and $K$ be an algebraic
number field such that $K \cap \Q(\zeta_{\alpha^l}) = \Q$. Assume
that
\begin{equation}
  \sum_{j \in X} c_j \cdot \zeta_{\alpha^l}^{-Nj^2+2aj} = 0
\end{equation}
where $X$ is a finite subset of $\Z$, $c_j \in K$ and $\alpha
\nmid a$. Then we have
\begin{equation}
  \sum_{j \in X : \;\; j \equiv 0 \mod \alpha^{l-1}}
  c_j \cdot \zeta_{\alpha^l}^{2aj} = 0
\end{equation}
if $\alpha \mid N$, or otherwise,
\begin{equation}
  \sum_{j \in X : \;\; Nj \equiv a \mod \alpha^{[(l+1)/2]}}
  c_j = 0.
\end{equation}
\end{lem}

\begin{proof}
If $\alpha \mid N$, taking the trace function of
$K(\zeta_{\alpha^l}) / K(\zeta_{\alpha})$ on both sides of the
equality assumed, we find from Lemma \ref{lem:A1} that
\begin{equation}
  \sum_{j \in X : \;\; j \equiv 0 \mod \alpha^{l-1}}
  c_j \cdot \zeta_{\alpha^l}^{2aj} = 0.
\end{equation}
Otherwise, choose $b \in \Z$ such that $bN \equiv 1 \mod
\alpha^l$. From the assumption we have
\begin{equation}
  \sum_{j\in X} c_j \cdot
  \zeta_{\alpha^l}^{-b(Nj-a)^2}
  = \zeta_{\alpha^l}^{-ba^2}
  \sum_{j\in X} c_j \cdot
  \zeta_{\alpha^l}^{-Nj^2+2aj}
  = 0.
\end{equation}
Taking the trace function of $K(\zeta_{\alpha^l}) /
K(\zeta_{\alpha})$ on both sides of the equality, we get
\begin{equation}
  \sum_{k=0}^{\alpha-1} \bigg(
  \sum_{j\in X : \;\; -b(Nj-a)^2 \equiv k\alpha^{l-1} \mod \alpha^l}
  c_j \bigg) \cdot \zeta_{\alpha}^{k}
  = 0.
\end{equation}
Since $\alpha$ is an odd prime, the congruence equation $x^2
\equiv -kN\alpha^{l-1} \mod \alpha^l$ has no solution for some
$0<k<\alpha$. It follows from Lemma \ref{lem:A2} that the
coefficient of $\zeta_\alpha^{k}$ in above summation identically
vanishes. In particular,
\begin{equation}
  \sum_{j\in X : \;\; -b(Nj-a)^2 \equiv 0 \mod \alpha^l} c_j = 0.
\end{equation}
Hence the lemma follows.
\end{proof}

\begin{proof}[Proof of Proposition]
Assume that $A_{p,q}^+(N,1)=0$ for
\begin{equation}
  p = \alpha_1^{l_1} \alpha_2^{l_2} \cdots \alpha_r^{l_r} \cdot
  \beta_1^{k_1} \beta_2^{k_2} \cdots \beta_s^{k_s},
\end{equation}
where the $\alpha_i$'s and $\beta_i$'s are distinct odd primes not
dividing and dividing $N$, respectively. Rewrite
$A_{p,q}^+(N,1)=0$ as
\begin{equation}
  -\frac14 \sum_{-pq<j<pq} j^2 \cdot
  \zeta_{4pq}^{-Nj^2+2qj} \cdot (\zeta_{2q}^{j} - \zeta_{2q}^{-j})
  = 0
\end{equation}
and choose $\sigma \in \Gal(\Q(\zeta_{8p})/\Q)$ so that
\begin{equation}
  \sigma(\zeta_{4pq})
  = \zeta_{4q} \cdot \zeta_{\alpha_1^{l_1}} \cdots \zeta_{\alpha_r^{l_r}} \cdot
  \zeta_{\beta_1^{k_1}} \cdots \zeta_{\beta_s^{k_s}}.
\end{equation}
Under the Galois action of $\sigma$, the equality becomes
\begin{equation}
  -\frac14 \sum_{-pq<j<pq} j^2 \cdot
  (\zeta_{4q} \cdot \zeta_{\alpha_1^{l_1}} \cdots \zeta_{\alpha_r^{l_r}} \cdot
  \zeta_{\beta_1^{k_1}} \cdots \zeta_{\beta_s^{k_s}})^{-Nj^2+2qj}
  \cdot (\zeta_{2q}^{pj} - \zeta_{2q}^{-pj})
  = 0.
\end{equation}
Put
\begin{equation}
  \alpha = \alpha_1^{[(l_1+1)/2]} \cdots \alpha_r^{[(l_r+1)/2]}, \quad
  \beta = \beta_1^{k_1-1} \cdots \beta_s^{k_s-1}, \quad
  p' = \beta_1 \cdots \beta_s.
\end{equation}
It follows from Lemma \ref{lem:A4} that
\begin{equation}\label{eqn:tmp1}
  -\frac14 \sum_{j \in X} j^2 \cdot
  (\zeta_{\beta_1^{k_1}} \cdots \zeta_{\beta_s^{k_s}})^{2qj}
  \cdot \zeta_{4q}^{-Nj^2+2qj} \cdot (\zeta_{2q}^{pj} - \zeta_{2q}^{-pj})
  = 0
\end{equation}
where
\begin{equation}
  X = \{ -pq<j<pq \mid
  Nj \equiv q \mod \alpha, \; j \equiv 0 \mod \beta \}.
\end{equation}
Now we apply the condition $q=2$. Notice that
\begin{equation}
  \zeta_8^{-Nj^2+4j} \cdot (\zeta_4^{pj} - \zeta_4^{-pj})
  = \left\{ \begin{array}{ll}
    0, & 2 \mid j, \\
    -\zeta_8^{-N} \cdot (-1)^{(j-1)/2} (\zeta_4^{p} - \zeta_4^{-p}), & 2 \nmid j.
  \end{array} \right.
\end{equation}
Dropping a nonzero factor, (\ref{eqn:tmp1}) becomes
\begin{equation}
  \sum_{j \in X : 2 \nmid j}
  j^2 \cdot \zeta_{p'}^{4j/\beta} \cdot (-1)^{(j-1)/2}
  = 0.
\end{equation}
Choose $0 \leq j_0 < 2\alpha$ so that $N\beta j_0 \equiv 2 \mod
\alpha$ and $j_0 \equiv 1 \mod 2$. Then the left hand side of
above equality, up to a sign, is
\begin{equation}
  \sum_{-p/\alpha\beta \leq j < p/\alpha\beta}
  \beta^2 (2\alpha j+j_0)^2 \cdot \zeta_{p'}^{8\alpha j+4j_0} \cdot (-1)^j
  = \frac {8p\beta \zeta_{p'}^{4j_0}} {1+\zeta_{p'}^{8\alpha}}
  \bigg( j_0 - \frac {2\alpha} {1+\zeta_{p'}^{-8\alpha}} \bigg).
\end{equation}
Therefore, we must have $p'=1$ and $j_0=\alpha$. But from the
choice of $j_0$, it follows that $\alpha=1$. Hence $p=1$, a
contradiction.
\end{proof}

\section{Proof of Theorem \ref{thm:wd}} \label{sec:wd}

Let $K$ denote the $r$-twisted Whitehead double of the torus knot
$T(p,q)$. Then
\begin{equation}
  J_{K,N}
  = \frac {1} {t^{N/2}-t^{-N/2}}
  \sum_{n=0}^{N-1} t^{rn(n+1)} \xi_{N,n} \hat{J}_{T(p,q),2n+1}
\end{equation}
where
\begin{equation}
  \hat{J}_{T(p,q),n} = (t^{n/2}-t^{-n/2}) J_{T(p,q),n}.
\end{equation}
Setting $t=\rt$, one notices that the denominator
$t^{N/2}-t^{-N/2}$ vanishes. Therefore, one has to apply
L'Hospital's rule, i.e. take derivative of both the denominator
and the numerator. It follows that
\begin{equation}\label{eqn:wd}
  J_{K,N}
  = \frac {-t^{-1/2}} {-N}
  \sum_{n=0}^{N-1} a_n^{4r-1} \sum_{i=0}^{N-1-n} S_{n,i}
  \Big( b_{n,i} \hat{J}_{T(p,q),2n+1} +
  t\frac{d}{dt} \hat{J}_{T(p,q),2n+1} \Big)
\end{equation}
where
\begin{equation}
  b_{n,i} = rn(n+1) - N(i+n) +
  \sum_{j=1}^{n} \Big (- \frac{N-i-j}{1-t^{-N+i+j}} - \frac{i+j}{1-t^{-i-j}}
  + \frac{j}{1-t^{-j}} \Big)
\end{equation}
and $a_n, S_{n,i}$ are same as Section \ref{sec:wl}.

Below, we follow the approach used in \cite{KT} to derive an
estimation of $\hat{J}_{T(p,q),n}$ and $t\frac{d}{dt}
\hat{J}_{T(p,q),n}$ in the form of (\ref{eqn:torus_estN}). For any
complex number $h$ with $\im(h)>0$, one has the integral formula
\begin{equation}\label{eqn:torus_int}
  \hat{J}_{T(p,q),n} (e^h)
  = e^{ -pq(n^2-1)\frac{h}{4} } \Big( \frac{pq}{\pi h} \Big)^{1/2}
  e^{ -(\frac{p}{q}+\frac{q}{p}) \frac{h}{4}}
  \int_C dz e^{pq \big( nz-\frac{z^2}{h} \big)} \tau(z)
\end{equation}
where the contour $C$ is given by the line $e^{\frac\pii4}\R$ and
\begin{equation}
  \tau(z) = \frac{(e^{pz}-e^{-pz})(e^{qz}-e^{-qz})}{e^{pqz}-e^{-pqz}}.
\end{equation}

\begin{lem}\label{lem:int_est}
For $h=\frac{2\pii}{N}$ we have
\begin{equation}\label{eqn:int_est}
  \frac{d^k}{dh^k} \int_C dz
  e^{pq \big( nz-\frac{z^2}{h} \big)} \tau(z)
  = -4\pi\sqrt{-1} \frac{1}{pq} \Big( \frac{N^2}{4pq} \Big)^k
  A_{p,q}^{(-1)^{n-1}}(N,k) + O(N^{2k-1/2})
\end{equation}
uniformly on $|n-N| < \frac{N}{2pq}$.
\end{lem}

\begin{proof}
Put $z_0 = \frac{n}{2}h = \frac{n}{N}\pii$. We have
\begin{equation}
  \quad \int_{C+z_0} dz
  e^{pq \big( nz-\frac{z^2}{h} \big)} z^{2k} \tau(z)
  = e^{pq\frac{z_0^2}{h}} \int_C dz
  e^{-pq\frac{z^2}{h}}
  (z+z_0)^{2k} \tau(z+z_0)
  = O(N^{-1/2})
\end{equation}
uniformly on $|n-N| < \frac{N}{2pq}$, since the function
$z^{2k}\tau(z)$ is bounded on the region
\begin{equation}
  z \in \Big\{ e^{\frac\pii4}x + y\pii \; \Big| \;
  x,y \in \R, |y-1|<\frac{1}{2pq} \Big\}.
\end{equation}
Counting the residues of the integrand at $\frac{j\pii}{pq}$,
$0<j<pq$, we also have
\begin{equation}
  (\int_C - \int_{C+z_0}) dz
  e^{pq \big( nz-\frac{z^2}{h} \big)} z^{2k} \tau(z)
  = -4 \Big( \frac{\pii}{pq} \Big)^{2k+1} A_{p,q}^{(-1)^{n-1}}(N,k).
\end{equation}
Therefore,
\begin{equation}
\begin{split}
  & \frac{d^k}{dh^k} \int_C dz
  e^{pq \big( nz-\frac{z^2}{h} \big)} \tau(z)
  = \Big( \frac{pq}{h^2} \Big)^k \int_C dz
  e^{pq \big( nz-\frac{z^2}{h} \big)} z^{2k} \tau(z) \\
  & \quad\quad\quad\quad\quad\quad
  = -4\pi\sqrt{-1} \frac{1}{pq} \Big( \frac{N^2}{4pq} \Big)^k
  A_{p,q}^{(-1)^{n-1}}(N,k)
  + O(N^{2k-1/2})
\end{split}
\end{equation}
uniformly on $|n-N| < \frac{N}{2pq}$.
\end{proof}

\begin{lem}\label{lem:torus_est}
For $t=\rt$ we have
\begin{equation}
  \hat{J}_{T(p,q),n} = O(1)
\end{equation}
and
\begin{equation}\label{eqn:torus_est2}
  t\frac{d}{dt} \hat{J}_{T(p,q),n}
  = -2 e^{ -pq\frac{n^2-1}{2N}\pii }  \frac{N^{5/2}}{(2pq)^{3/2}}
  e^{ -(\frac{p}{q}+\frac{q}{p}) \frac{\pii}{2N} + \frac\pii4 }
  A_{p,q}^{(-1)^{n-1}}(N,1) + O(N^2)
\end{equation}
uniformly on $|n-N| < \frac{N}{2pq}$.
\end{lem}

\begin{proof}
From (\ref{eqn:torus_int}) and Lemma \ref{lem:int_est} we have
\begin{equation}\label{eqn:torus_est1}
  \hat{J}_{T(p,q),n}
  = -4 e^{ -pq\frac{n^2-1}{2N}\pii } \frac{N^{1/2}}{(2pq)^{1/2}}
  e^{ -(\frac{p}{q}+\frac{q}{p}) \frac{\pii}{2N} + \frac\pii4 }
  A_{p,q}^{(-1)^{n-1}}(N,0) + O(1).
\end{equation}
It follows from the periodicity of $A_{p,q}^\pm$ and the identity
\begin{equation}
  \hJ{T(p,q)}{N} = 0 \cdot \J{T(p,q)} = 0
\end{equation}
that $A_{p,q}^{(-1)^{N-1}}(N,0)$, hence by (\ref{eqn:A})
$A_{p,q}^\pm(N,0)$, identically vanishes. So, the leading term of
the right hand side of (\ref{eqn:torus_est1}) is zero. Then,
applying Lemma \ref{lem:int_est} to the derivative of
(\ref{eqn:torus_int}), one obtains the second estimation.
\end{proof}

Now we conclude the proof of the theorem. It is clear that
\begin{equation}
  b_{n,i} \hat{J}_{T(p,q),2n+1} + t\frac{d}{dt} \hat{J}_{T(p,q),2n+1}
  = O(N^3) O(N) + O(N^3) = O(N^4)
\end{equation}
uniformly on $0 \leq n,i,n+i < N$. Moreover, for $0<\alpha<1$,
since the function
\begin{equation}
  \frac{x}{1-e^{-2\sqrt{-1}x}}
  = \frac{e^{\sqrt{-1}x}x}{2\sqrt{-1}\sin x}
\end{equation}
is bounded on $x \in [0,\alpha\pi]$, we have
\begin{equation}
  b_{n,i} = O(N^2)
\end{equation}
uniformly on $0 < n,i,n+i < \alpha N$. It follows from Lemma
\ref{lem:torus_est} that
\begin{equation}\label{eqn:torus_est1a}
  b_{n,i}\hat{J}_{T(p,q),2n+1} = O(N^2)
\end{equation}
and
\begin{equation}\label{eqn:torus_est1b}
  t\frac{d}{dt} \hat{J}_{T(p,q),2n+1}
  = -2 (a_n)^{-4pq} \frac{N^{5/2}}{(2pq)^{3/2}}
  e^{ -(\frac{p}{q}+\frac{q}{p}) \frac{\pii}{2N} + \frac\pii4 }
  A_{p,q}^+(N,1) + O(N^2)
\end{equation}
uniformly on $|n-\hN|+|i-\qN| < \frac{N}{4pq}$.

Therefore, in the case that $q=2$, by Lemma \ref{lem:main1}, Lemma
\ref{lem:main2} and Proposition \ref{prop:novanish}, in the same
notations as the lemmas we have
\begin{equation}
\begin{split}
  & \sum_{|n-\hN|+|i-\qN| \geq N^\delta} a_n^{4r-1} S_{n,i}
  \Big( b_{n,i} \hat{J}_{T(p,q),2n+1} + t\frac{d}{dt} \hat{J}_{T(p,q),2n+1} \Big)
  = N^6 e^{-\epsilon N^{2\delta-1}} S_{[\hN],[\qN]} O(1), \\
  & \sum_{|n-\hN|+|i-\qN| < N^\delta} a_n^{4r-1} S_{n,i}
  \Big( b_{n,i} \hat{J}_{T(p,q),2n+1} + t\frac{d}{dt} \hat{J}_{T(p,q),2n+1} \Big)
  = N^{7/2} S_{[\hN],[\qN]} e^{O(1)},
\end{split}
\end{equation}
hence
\begin{equation}
\begin{split}\label{eqn:Jwd_est}
  2\pi \logJ{K}
  & = 2\pi \log (N^{5/2} S_{[\hN],[\qN]}) + O(1) \\
  & = 8 L(\frac{\pi}{4}) \cdot N + 4\pi \log N + O(1) \\
  & = v_3\sv{K} \cdot N + 4\pi \log N + O(1)
\end{split}
\end{equation}
as $N\to\infty$.

\section {Concluding remarks}

Although the proof of Theorem \ref{thm:wd} depends on the simple
nature of the Whitehead doubles of torus knots, the approach still
works for the satellite knots to which the colored Jones
polynomials of the associated companion knot and pattern link
satisfy certain mild conditions. Meanwhile, one has to deal with
several problems.

For example, as shown in expression (\ref{eqn:Jwd}), although the
volume conjecture itself is only concerned with the value of the
colored Jones polynomial $J_{K,N}$ at the $N$-th root of unity,
the values at other roots of unity become crucial once the
satellite knots of $K$ are involved.

A more challenging problem is due to the estimations
(\ref{eqn:torus_est1a}) and (\ref{eqn:torus_est1b}), which have
enabled us to neglect the term $b_{n,i}\hat{J}_{T(p,q),2n+1}$ in
summation (\ref{eqn:wd}). Note that the derivative of the
polynomial
\begin{equation}
  \hat{J}_{K,N} = (t^{N/2}-t^{-N/2}) J_{K,N}
\end{equation}
is related to $J_{K,N}$ via the identity
\begin{equation}
  t\frac{d}{dt}\hJ{K}{N} = -N \J{K}.
\end{equation}
Therefore, the term $t\frac{d}{dt} \hat{J}_{T(p,q),2n+1}$ in
summation (\ref{eqn:wd}), in fact, plays the role of $\J{T(p,q)}$.
Hence it is quite natural to see the term
$b_{n,i}\hat{J}_{T(p,q),2n+1}$ be suppressed. Following this
observation, when the Whitehead doubles of general knots are
considered, it is reasonable to expect that a similar suppression
also happens.

\begin{conj}
For every nontrivial knot $K$ we have
\begin{equation}
  \frac {\hJ{K}{2n+1}}
  {t\frac{d}{dt} \hJ{K}{2n+1}}
  = o(N^{-2})
\end{equation}
uniformly on $|n-\hN| < N^\delta$ for some $\frac12 < \delta <
\frac23$.
\end{conj}

Note that the conjecture excludes the case of unknot, for which
the statement of the conjecture is obviously false. Indeed, the
Whitehead doubles of unknot are no longer satellite knots but the
so called twist knots (including unknot, trefoil, figure 8, etc.),
whose complements always admit a volume strictly smaller than that
of Whitehead link. In the sequel, the conjecture has an
interesting implication: the colored Jones polynomial detects
unknot.

\section*{Acknowledgement}

The author is grateful to Xiao-Song Lin for enlightening
discussions. He also thanks Lingquan Kong who pointed out to the
author the approach to proving the nonvanishing proposition.


\begin{thebibliography}{999}
\bibitem{Hikami} K. Hikami,
    Volume conjecture and asymptotic expansion of $q$-series,
    Exper. Math. 12 (2003), 319--337.
\bibitem{Kashaev} R. M. Kashaev,
    The hyperbolic volume of knots from quantum dilogarithm,
    Lett. Math. Phys. 39 (1997), 269--275.
\bibitem{KT} R. M. Kashaev, O. Tirkkonen,
    Proof of the volume conjecture for torus knots,
    J. Math. Sci. 115 (2003), 2033--2036.
\bibitem{Morton} H. R. Morton,
    The coloured Jones function and Alexander polynomial for torus knots,
    Proc. Cambr. Phil. Soc. 117 (1995), 129--135.
\bibitem{MM} H. Murakami, J. Murakami,
    The colored Jones polynomials and the simplicial volume of a knot,
    Acta Math. 186 (2001), 85--104.
\bibitem{MMOTY} H. Murakami, J. Murakami, M. Okamoto, T. Takata, Y. Yokota,
    Kashaev's conjecture and the Chern-Simons invariants of knots and links,
    Exper. Math. 11 (2002), 427--435
\bibitem{Thurston} W. P. Thurston,
    Three-dimensional manifolds, Kleinian groups, and hyperbolic geometry,
    Bull. A.M.S. 6 (1982), 357--381.
\end{thebibliography}
\end{document}